\documentclass{amsart}

\usepackage{amsmath,amscd,amssymb, epsfig,psfrag}

\newtheorem{theorem}{Theorem}[section]

\theoremstyle{definition}

\newtheorem{conjecture}[theorem]{Conjecture}

\newcommand{\ZZ}{\mathbb{Z}}

\newcommand{\QQ}{\mathbb{Q}}

\newcommand{\RR}{\mathbb{R}}

\def\d{\partial}

\def\e{\epsilon}
\def\g{\gamma}

\def\G{\Gamma}

\def\S{\Sigma}

\def\sech{\mbox{\rm{sech}}}
\def\arctanh{\mbox{\rm{arctanh}}}

\def\Isom{\mbox{Isom}}

\def\O{\mbox{\rm{O}}}
\def\PO{\mbox{\rm{PO}}}
\def\OO{\mathcal{O}}

\def\HH{\mathbb{H}}

\edef\t@mp{\catcode`\noexpand\~=\the\catcode`\~}%
    \catcode`\~=12
    \def\tild@{~}%
    \t@mp
\begin{document}
\Large

\title{\Large Systoles of Hyperbolic 4-manifolds} 

\author{Ian Agol}
\thanks{Agol partially supported by NSF grant DMS-0504975 and the Guggenheim Foundation}

\email{agol@math.uic.edu}
\address{MSCS UIC 322 SEO, m/c 249\\
        851 S. Morgan St.\\
        Chicago, IL 60607-7045}

%

\date{%
 May 18, 2006}


\begin{abstract} 
\Large We prove that for any $\e>0$, there exists a closed hyperbolic 4-manifold with a closed geodesic of length $< \e$. 
\end{abstract} 



\maketitle

\section{Introduction}
It has been known for a long time that closed hyperbolic surfaces
and 3-manifolds may have arbitarily short geodesics. 
This follows in 2-dimensions by an explicit construction,
and holds in 3 dimensions by Thurston's hyperbolic 
Dehn surgery theorem \cite[Thm. 5.8.2]{Th}. In this note, we prove the 
existence of closed hyperbolic 4-manifolds which have
arbitrarily short geodesics. It is conjectured that there
is a uniform lower bound on the length of a systole in 
arithmetic hyperbolic manifolds. This would follow from
Lehmer's conjecture \cite[\textsection10]{Gelander04}. Examples
of non-arithmetic hyperbolic manifolds in higher dimensions  
come from the  method of ``inter-breeding'' introduced 
by Gromov and Piatetski-Shapiro \cite{GromovPS88}, by
taking hyperbolic manifolds with geodesic boundary produced using
arithmetic methods, and gluing them to obtain a  non-arithmetic
manifold. Our result makes use of a variation on their method, which might
best be described as ``inbreeding''. The method
would extend to all dimensions if a conjecture about 
the fundamental groups of arithmetic hyperbolic manifolds were known.
We speculate on this conjecture in the last section.

\section{Subgroup separability}
Let $G$ be an infinite group, and $A<G$ a subgroup. We say that $A$ is
{\it separable} in $G$ if $A$ is the intersection of all finite index subgroups
of $G$ which contain $A$, that is
\begin{equation}
A= \underset{A\leq B\leq G,  [G:B]<\infty}{\cap} B.
\end{equation} 
We say that a discrete group $\G<\Isom(\HH^{n})$ is  {\it GFERF} (short for {\it Geometrically
Finite Extended Residually Finite}) if every geometrically finite subgroup $A<\G$
is separable in $\G$. More generally, a group $G$ is {\it LERF} if every 
finitely generated subgroup $A<G$ is separable in $G$. If $\G$ is LERF,
then since geometrically finite subgroups are finitely generated, this would
imply that $\G$ is GFERF. Unfortunately, the converse is not necessarily
true, since there may be finitely generated subgroups of $\G$ which are 
not geometrically finite. 
Scott showed that $A< \G$ is separable if and only if for any 
compact subset $C\subset \HH^{4}/A$, there exists $\G_{1}< \G$, $[\G:\G_{1}]< \infty$, and 
$C \hookrightarrow \HH^{4}/\G_{1}$ embeds under the covering map.
We will be using a fact due to Scott \cite{S}
that the group generated by reflections in the right-angled 120-cell in $\HH^{4}$
is GFERF; for a proof see \cite{ALR}. 

\section{Systoles}
\begin{theorem} \label{systole}
There exist closed hyperbolic 4-manifolds with arbitrarily short geodesics. 

\end{theorem}

\begin{proof}
The examples will come from cutting and pasting certain covers of a 4-dimensional Coxeter orbifold.
The right-angled 120-cell $D\subset \HH^{4}$ is the fundamental domain for a reflection group. 
Let $\G_{D}$ be the group generated by reflections
in the faces of $D$. Let $\OO$ be the ring of integers in $\QQ(\sqrt{5})$.
Then $\G_{D}$ is commesurable with $\PO(f;\OO)$, where $f$ is the 5-dimensional quadratic
form $\langle 1,1,1,1,-\phi\rangle$, where $\phi=\frac{1+\sqrt{5}}{2}$ \cite[Lemma 3.3]{ALR}. Let $P \subset \HH^{4}$ be
a 3-dimensional geodesic subspace, such that $H=\Isom(P)\cap \G_{D}$ is a cocompact subgroup
of $\Isom(P)$ (where we embed $\Isom(P)<\Isom(\HH^{4})$ in the
natural fashion). If we identify $\HH^{4}$ with a component of the hyperboloid $f(x)=-1$, then we
find such a $P$ by letting $P=\HH^{4}\cap v^{\perp}$ (with respect to the inner product defined
by $f$) where $v \in \QQ(\sqrt{5})^{5}$, $f(v)>0$. Let $\G$ be a finite index torsion-free subgroup of $\G_{D}$, which exists by 
Selberg's lemma.  
Now, $Comm(\G) > \PO(f, \QQ(\sqrt{5}))$, so $Comm(\G)$ is dense in $\PO(f,\RR)=\Isom(\HH^{4})$.   
Thus, for any $\e>0$,  we may find $\g \in Comm(\G)$ such that $\g(P)\cap P = \emptyset$, and $d(P,\g(P)) < \frac{\e}{2}$.
The plane $\g(P)$ is stabilized by $(\g H \g^{-1})\cap \G$, which is cocompact in
$\Isom(\g(P))$, since $ [\G: (\g \G \g^{-1} )\cap \G]<\infty$. Thus, $H_{\g}=\Isom(\g(P)) \cap \G$ is cocompact
in $\Isom(\g(P))$, since $(\g H \g^{-1})\cap\G < H_{\g}$.  
Let $g\subset \HH^{4}$ be a geodesic segment perpendicular to $P$ and $\g(P)$, and with endpoints
$p_{1}=g\cap P$, and $p_{2}= g\cap \g(P)$. Let $\rho: \RR\to \RR$ be the function $\rho(x)= \arctanh\ \sech x$.
Using residual finiteness, choose $H_{1}< H$ such that $d(p_{1},h(p_{1})) > 2\rho(l(g)/2)$, for all $h \in H_{1}-\{1\}$. 
Similarly, choose $H_{2}< H_{\g}$ such that $d(p_{2}, h(p_{2}))> 2\rho(l(g)/2)$, for all $h\in H_{2}-\{1\}$.
Let $\S_{1}= P/H_{1}$, and $\S_{2}=\g(P)/H_{2}$. 

Claim: $G=\langle H_{1}, H_{2}\rangle \cong H_{1} \ast H_{2}$. Moreover, $\HH^{4}/G$ is geometrically
finite. 

To prove the claim, let $E_{i}\subset \HH^{4}$ be a Dirichlet domain about $p_{i}$ with respect to the group $H_{i}$. 
Let $L$ be the 3-plane which is
the perpendicular bisector of $g$. Let $pr_{1}: \HH^{4} \to P$, and $pr_{2}: \HH^{4} \to \g(P)$. 
Then $pr_{i}(L)$ is a disk about $p_{i}$ of radius $\rho(l(g)/2)$ in 
$P$ or $\g(P)$ (see Figure \ref{projection}). 
\begin{figure}[htb] 
	\begin{center}
	\psfrag{a}{$\rho(l/2)$}
	\psfrag{b}{$l/2$}
	\psfrag{c}{$g$}
	\psfrag{d}{$p_{1}$}
	\psfrag{e}{$p_{2}$}
	\psfrag{f}{$P$}
	\psfrag{g}{$\g(P)$}
	\psfrag{h}{$L$}
	\psfrag{i}{$E_{1}$}
	\psfrag{j}{$E_{2}$}
	\epsfig{file=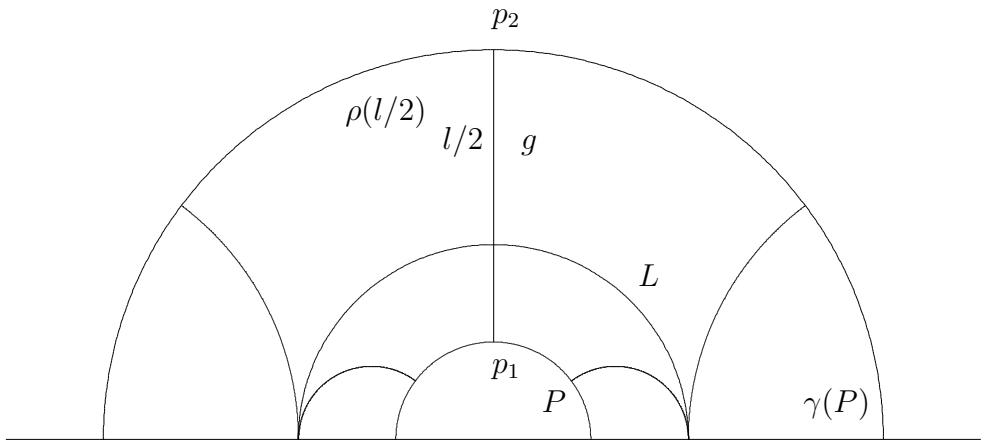,width=\textwidth,height=.5\textwidth}
	\caption{\label{projection} Combining fundamental domains}
	\end{center}
\end{figure}

Since $E_{i}$ is a Dirichlet domain, it must contain $pr_{i}(L)$, and
therefore $E_{i}$ contains $L$. Thus,  $\d E_{1} \cap \d E_{2}=\emptyset$ since they are separated by the hyperplane $L$. Thus, $E_{1}\cap E_{2}$
will be  a finite-sided fundamental domain for $G$, and thus $G$ is geometrically finite (see \cite{Bowditch93} for
various equivalent notions of geometric finiteness). 
Topologically, $\HH^{4}/G = (\HH^{4}/H_{1}) \#_{L} (\HH^{4}/H_{2})$, so $ G \cong H_{1}\ast H_{2}$.   
 
Let $U= \S_{1} \cup_{p_{1}} g \cup_{p_{2}} \S_{2}$. Then $U$ is an embedded compact spine of $\HH^{4}/G$.  
Now, we use the fact from theorem  3.1 \cite{ALR}, that $G$ is a separable subgroup of $\G$. 
By Scott's separability criterion, we see that we
may embed $U$ in $\HH^{4}/\G_{1}$, for some finite index subgroup $\G_{1}<\G$. Thus, we have $\S_{1}\cup \S_{2} \subset \HH^{4}/\G_{1}$.
Let $N=(\HH^{4}/\G_{1} )\backslash (\S_{1} \cup \S_{2})$, and let $M=DN$, the double of $N$
along its boundary. Since $g\subset N$ is a geodesic arc orthogonal to $\d N$, we have the double
of $g$  $D(g) \subset M$ is a closed geodesic in $M$ of length $<\e$. 
\end{proof}

\section{Conclusion}
\begin{conjecture}
There exists closed hyperbolic $n$-manifolds with arbitrarily short geodesics. 
\end{conjecture}
Hyperbolic lattices that are subgroup separable on geometrically finite subgroups are 
called {\it GFERF}, short for {\it Geometrically Finite Extended Residual Finite}.  
This conjecture would follow from the following conjecture, by the same proof as the main theorem.
\begin{conjecture}
There exist compact arithmetic hyperbolic $n$-manifolds which are defined by a quadratic form, and 
which are GFERF. 
\end{conjecture}
Unfortunately, there does not exist a compact right-angled polyhedron in $\HH^{n}$, $n\geq 5$, so
the strategy of proof in \cite{ALR, S} will not work in general. 
By the remark after \cite[Lemma 3.4]{ALR}, we know that $\O(8,1;\ZZ)$ is GFERF. 
The above conjecture would hold if we knew $\O(n,1;\ZZ)$ were GFERF for all $n$, since one
may embed (up to finite index) any cocompact arithmetic lattice defined by a quadratic form into $\O(n,1;\ZZ)$
for some $n$ by the restriction of scalars and stabilization. The following
theorem is proven the same as Theorem \ref{systole}:

\begin{theorem}
There exist finite volume hyperbolic $n$-manifolds with arbitrarily short geodesics for $n\leq 8$. 
\end{theorem}

At most finitely many of the manifolds produced using Theorem \ref{systole} will be
arithmetic. This follows because the groups will lie in $\O(f, \QQ(\sqrt{5}))$. 
The integral real eigenvalues of the matrices in $\O(f,\QQ(\sqrt{5}))$ will
be bounded away from 1, since they have an integral minimal polynomial of degree
at most 10. Thus, the length of a geodesic of an arithmetic subgroup of $\O(f,\QQ(\sqrt{5}))$
will be bounded away from 0, which implies that at most finitely many examples 
from Theorem \ref{systole} may be arithmetic. This method for proving the existence
of non-arithmetic uniform lattices is slightly different than the method of \cite{GromovPS88}, since
instead of breeding subgroups of incommensurable arithmetic lattices, it breeds a subgroup
of an arithmetic lattice with itself. It's possible that this ``inbreeding'' method could produce
non-arithmetic lattices in any dimension.  

\bibliographystyle{hamsplain}
\bibliography{4Dhyperbolicsystole}

\def\cprime{$'$} \def\cprime{$'$} \def\cprime{$'$}
\providecommand{\bysame}{\leavevmode\hbox to3em{\hrulefill}\thinspace}
\providecommand{\href}[2]{#2}
\begin{thebibliography}{1}

\bibitem{ALR}
I.~Agol, D.~D. Long, and A.~W. Reid, \emph{The {B}ianchi groups are separable
  on geometrically finite subgroups}, Ann. of Math. (2) \textbf{153} (2001),
  no.~3, 599--621.

\bibitem{Bowditch93}
B.~H. Bowditch, \emph{Geometrical finiteness for hyperbolic groups}, J. Funct.
  Anal. \textbf{113} (1993), no.~2, 245--317.

\bibitem{Gelander04}
Tsachik Gelander, \emph{Homotopy type and volume of locally symmetric
  manifolds}, Duke Math. J. \textbf{124} (2004), no.~3, 459--515.

\bibitem{GromovPS88}
M.~Gromov and I.~Piatetski-Shapiro, \emph{Nonarithmetic groups in {L}obachevsky
  spaces}, Inst. Hautes \'Etudes Sci. Publ. Math. (1988), no.~66, 93--103.

\bibitem{S}
Peter Scott, \emph{Subgroups of surface groups are almost geometric}, Journal
  of the London Mathematical Society. Second Series \textbf{17} (1978), no.~3,
  355--565.

\bibitem{Th}
William~P. Thurston, \emph{The geometry and topology of 3-manifolds}, Lecture
  notes from Princeton University, 1978--80.

\end{thebibliography}

\end{document}